\documentstyle{article}
\begin{document}

\baselineskip 15pt
\parindent=1em
\hsize=12.3 cm \textwidth=12.3 cm
\vsize=18.5 cm \textwidth=18.5 cm

\def\supt{{\rm supt}}
\def\dom{{\rm dom}}
\def\bfone{{\bf 1}}
\def\Gen{{\rm Gen}}
\def\rk{{\rm rk}}
\def\cali{{\cal I}}
\def\calj{{\cal J}}
\def\Spec{{\rm Spec}}
\def\cf{{\rm cf}}

\title{Understanding preservation theorems,  II}
\author{
Chaz Schlindwein \\
Department of Mathematics and Computing \\
Lander University \\
Greenwood, South Carolina 29649 USA\\
{\tt cschlind@lander.edu}}

\maketitle

\def\forces{\mathbin{\parallel\mkern-9mu-}}
\def\notforces{\,\nobreak\not\nobreak\!\nobreak\forces}

\def\restr{\,\hbox{\vrule height8pt width.4pt depth0pt
   \vrule height7.75pt width0.3pt depth-7.5pt\hskip-.2pt
   \vrule height7.5pt width0.3pt depth-7.25pt\hskip-.2pt
   \vrule height7.25pt width0.3pt depth-7pt\hskip-.2pt
   \vrule height7pt width0.3pt depth-6.75pt\hskip-.2pt
   \vrule height6.75pt width0.3pt depth-6.5pt\hskip-.2pt
   \vrule height6.5pt width0.3pt depth-6.25pt\hskip-.2pt
   \vrule height6.25pt width0.3pt depth-6pt\hskip-.2pt
   \vrule height6pt width0.3pt depth-5.75pt\hskip-.2pt
   \vrule height5.75pt width0.3pt depth-5.5pt\hskip-.2pt
   \vrule height5.5pt width0.3pt depth-5.25pt}\,}

   \centerline{{\bf Abstract}}

   We present an exposition of much of Sections VI.3 and XVIII.3 from Shelah's book {\em Proper and 
Improper Forcing}. This covers numerous preservation theorems for countable support iterations
of proper forcing, including preservation of the property ``no new random reals
over $V$,'' the property ``reals of the ground model form
 a non-meager set,'' 
 the property ``every dense open set contains a dense open set of the ground model,''
and preservation theorems related to the weak bounding property, the weak ${}^\omega\omega$-bounding
property, and the property ``the set of reals of the ground model has positive
outer measure.''

\eject

\section{Introduction}

This is the fourth of a sequence of papers giving an exposition of portions of 
Shelah's book, {\em Proper and Improper Forcing} [9].  
The earlier papers were [6], [7], and [8],
which cover sections 2 through 8 of [9, Chapter XI], sections 2 and 3 of [9, Chapter XV],
and sections 1 and 2 of [9, Chapter VI], respectively. 

In this paper, we give an exposition of much of [9, Sections VI.3 and XVIII.3] dealing with preservation
theorems. We include proofs of the preservation, under countable support
iteration of proper forcing, of the property
``no new random reals,'' the property ``every open dense set contains an old open dense set,'' 
 the property of non-meagerness of the reals
of the ground model, and preservation theorems related to
weak bounding, weak ${}^\omega\omega$-bounding, and ``the set of reals of the ground model
has positive outer measure.''
 
Another treatment of preservation theorems, using different methods, is given in [2], [3].
The results of [9, Section VI.3] included here as Theorem 2.5, Theorem 3.5, and Theorem 4.13 may also be derived as corollaries of [1, Theorem 6.1.18]; the proof there is 
essentially the same as the ones given by Shelah in [9, Section VI.3].

\section{Preservation of weak bounding}

The most important tool in the study of preservation theorems for countable support forcing iterations is the Proper Iteration Lemma.  Here, and throughout this paper, $P_{\alpha,\kappa}$ is characterized by

\medskip

\centerline{
$V[G_{P_\alpha}]\models``P_{\alpha,\kappa}=\{p\restr[\alpha,\kappa)\,\colon p\in P_\kappa$ and $p\restr \alpha\in G_{P_\alpha}\}$.''}

\proclaim Theorem 2.1 (Proper Iteration Lemma, Shelah). Suppose\/ $\langle P_\eta\,
\colon\allowbreak\eta\leq\kappa\rangle$
is a countable support forcing iteration based on\/
$\langle Q_\eta\,\colon\allowbreak\eta<\kappa\rangle$ and
for every\/ $\eta<\kappa$ we have that\/ ${\bf 1}\forces_{P_\eta}`` Q_\eta$ {\rm is
proper.''} Suppose also that\/ $\alpha<\kappa$ and\/
$\lambda$ is a sufficiently large regular cardinal and\/ $N$ is a
countable elementary submodel of\/ $H_\lambda$ and\/
$\{P_\kappa,\alpha\}\in N$ and\/
$p\in P_\alpha$ is\/ $N$-generic and\/ {\rm
$p\forces``q\in  P_{\alpha,\kappa}\cap
N[G_{P_\alpha}]$.''}
Then there is\/ $r\in P_\kappa$ such that\/
$r$ is\/ $N$-generic and\/ $r\restr\alpha=p$ and\/ {\rm
$p\forces``r\restr[\alpha,\kappa)\leq
q$.''}

Proof: See (e.g.) [8, Theorem 2.1].

We deal first with the weak bounding property.

\proclaim Definition 2.2. Suppose $A$ and $B$ are sets of integers.  We say $A\subseteq^* B $ iff $\{n\in A\,\colon\allowbreak n\notin B\}$ is finite.

\proclaim Definition 2.3. Suppose ${\cal P}\subseteq[\omega]^{\aleph_0}$ is a filter.
 We say ${\cal P}$ is a P-filter iff ${\cal P}$ contains all co-finite subsets of
$\omega$, and  $(\forall{\cal U}\in[{\cal P}]^{\aleph_0})\allowbreak(\exists A\in{\cal P})\allowbreak
(\forall B\in{\cal U})\allowbreak(A\subseteq^* B)$.

\proclaim Definition 2.4. Suppose ${\cal P}$ is a P-filter and\/ $P$ is a forcing notion. 
We say that $P$ is weakly ${\cal P}$-bounding iff\/ 
{\rm ${\bf 1}\forces_P``(\forall A\in[\omega]^{\aleph_0})
\allowbreak(\exists B\in {\cal P})\allowbreak(A\not\subseteq^* B)$.''}

The following Theorem is [9, Conclusion VI.3.17(1)].

\proclaim Theorem 2.5.  Suppose $\kappa$ is a limit ordinal
and ${\cal P}$ is a P-filter and\/
$\langle P_\eta\,\colon\eta\leq\kappa\rangle$ is a countable support forcing iteration based on
$\langle Q_\eta\colon\allowbreak\eta<\kappa\rangle$.  Suppose for every $\eta<\kappa$ 
we have\/ $P_\eta$ is weakly ${\cal P}$-bounding and\/ {\rm
${\bf 1}\forces_{P_\eta}``Q_\eta$ is proper.''}
Then\/ $P_\kappa$ is weakly ${\cal P}$-bounding.

Proof:
This is clear if $\kappa$ has uncountable cofinality, so assume
 $\cf(\kappa)=\omega$.

Suppose $p\in P_\kappa$ and $A$ is a $P_\kappa$-name and 
$p\forces``A\in[\omega]^{\aleph_0}$.'' Let $\lambda$ be a 
sufficiently large regular cardinal and
$N$ a countable elementary substructure of $H_\lambda$ such that 
$\{P_\kappa,{\cal P},A,p\}\in N$.

Let $\langle\alpha_k\,\colon\allowbreak k\in\omega\rangle\in N$ be an increasing sequence cofinal in $\kappa$
such that $\alpha_0=0$.
Fix $B\in{\cal P}$ such that $(\forall X\in{\cal P}\cap N)
\allowbreak(B\subseteq^* X)$.  It suffices to show $p\notforces``A\subseteq^* B$.''

Build $\langle q_k,p_k,m_k\,\colon k\in\omega\rangle$ such that $q_0=p$ and 
for every $k\in\omega$ we have that each of the following holds:

(1) $p_k\in P_{\alpha_k}$ is $N$-generic, and

(2) $p_k\forces``q_{k+1}\in P_{\alpha_k,\kappa}\cap N[G_{P_{\alpha_k}}]$ and $q_{k+1}\leq q_k\restr[\alpha_k,\kappa)$,'' and

(3) $m_k$ is a $P_{\alpha_k}$-name for an integer and $p_k\forces``$if $k>0$ then 
$m_k>m_{k-1}$,'' and

(4) $p_k\forces``q_{k+1}\forces`m_k\notin  B$ and $m_k\in A$,'\thinspace'' and

(5) $p_{k+1}\restr\alpha_k =p_k$, and

(6) $p_k\forces``p_{k+1}\restr[\alpha_k,\alpha_{k+1})\leq q_{k+1}\restr\alpha_{k+1}$.''

The construction proceeds as follows. Given $p_k$, $q_k$, and $m_{k-1}$, work in
 $V[G_{P_{\alpha_k}}]$ with $p_k\in G_{P_{\alpha_k}}$.

Build $A_k\in[\omega]^{\aleph_0}\cap N[G_{P_{\alpha_k}}]$ and $\langle q^i_k\,\colon\allowbreak i\in\omega\rangle\in N[G_{P_{\alpha_k}}]$ such that 
$q^0_k=q_{k}\restr[\alpha_k,\kappa)$, and for every $i\in\omega$ we have 
 $q^{i+1}_k\leq q^i_k$, and
 $q_k^{i+1}\forces``i\in A $ iff $i\in A_k$.''

Using the hypothesis on $P_{\alpha_k}$ we may choose $B_k\in{\cal P}$ 
such that $A_k\not\subseteq^* B_k$.  By elementarity we may assume
$B_k\in N[G_{P_{\alpha_k}}]$.  Because $p_k$ is $N$-generic, we have
$B_k\in V\cap N[G_{P_{\alpha_k}}]=N$.

Because $B_k\in N$ we have $B\subseteq^* B_k$, and hence $A_k\not\subseteq^* B$.  Therefore we may
 choose $m_k\in\omega$ such that  $m_k\in A_k$ and $m_k\notin B$ and if $k>0$ then $m_k>m_{k-1}$.

Let $q_{k+1}=q_k^{m_k+1}$. Clearly (2), (3), and (4) are satisfied.

Using the Proper Iteration Lemma we may choose $p_{k+1}$ satisfying (1), (5), and (6).

This completes the recursive construction.

Let $r\in P_\kappa$ be such that for every $k\in\omega$ we have $r\restr \alpha_k
=p_k$.

Suppose, towards a contradiction, that $r'\leq r$ and $r'\forces``A\subseteq^* B$.''  
Fix $P_\kappa$-names $n$ and $k$ such that
$r'\forces``A\subseteq n\cup B$, and $m_k>n$.''

By strengthening $r'$ we may assume that $k$ and $m_k$ are integers rather than merely  names.

Because $r'\leq(p_{k+1},q_{k+1})$ we have
$r'\forces``m_k\in A-n\subseteq B$, and $m_k\notin B$.''
This is a contradiction.

The Theorem is established.

\proclaim Lemma 2.8.  Suppose $P$ is weakly\/
 ${\cal P}$-bounding and\/ {\rm ${\bf 1}\forces_P``Q $ is almost ${\cal P}$-bounding.''} 
 Then $P*Q$ is weakly ${\cal P}$-bounding.

Proof: Suppose $(p,q)\forces_{P*Q}``A\in[\omega]^{\aleph_0}$.'' Take $q'$ and $B$ in $V^P$ such that $p\forces``B\in{\cal P}$ and $q'\leq q$ and

(*) $q'\forces`(\forall Y\in[\omega]^{\aleph_0}\cap V[G_P])(A\cap Y\not\subseteq^* B$).'\thinspace''

Take $p'\leq p$ and $B'\in{\cal P}$ such that $p'\forces``B\not\subseteq^* B'$.''
By (*) we have $(p',q')\forces``A\cap(B-B')\not\subseteq^* B$.'' Hence
$(p',q')\forces``A\not\subseteq^* B'$.''

The Lemma is established.

\proclaim Theorem 2.9.  Suppose ${\cal P}$ is a P-filter and\/ $\langle 
P_\eta\,\colon\allowbreak\eta\leq\kappa\rangle$ is
a countable support forcing iteration based on $\langle Q_\eta\,\colon\allowbreak\eta
<\kappa\rangle$.  Suppose for every $\eta<\kappa$ we have\/ {\rm 
${\bf 1}\forces_{P_\eta}``Q_\eta$
is proper and almost ${\cal P}$-bounding.''}
Then $P_\kappa$ is weakly ${\cal P}$-bounding.

Proof: By Theorem 2.5 and Lemma 2.8.

\section{Preservation of weakly ${}^\omega\omega$-bounding}

In this section we give an exposition of a preservation theorem, due to Shelah, concerning
the weak ${}^\omega\omega$-bounding property.

\proclaim Definition 3.1. Suppose $f$ and $g$ are in ${}^\omega\omega$.  We say $f\leq^* g$ iff $(\exists n\in\omega)\allowbreak(\forall k > n)\allowbreak( f(k)\leq g(k))$.

\proclaim Definition 3.2. Suppose $F\subseteq{}^\omega\omega$ and $g\in{}^\omega\omega$.  We say that $g$ bounds $F$ iff $(\forall f\in F)\allowbreak(f\leq^* g)$.

\proclaim Definition 3.3. Suppose $P$ is a forcing notion. We say that $P$ is weakly ${}^\omega\omega$-bounding iff\/ {\rm ${\bf 1}\forces``(\forall f\in{}^\omega\omega)\allowbreak(\exists g\in{}^\omega\omega\cap V)\allowbreak(g\not\leq^* f)$.''}

The following Theorem is [9, Conclusion VI.3.17(2)].

\proclaim Theorem 3.4.  Suppose $\kappa$ is a limit ordinal and $\langle P_\eta\,\colon\eta\leq\kappa\rangle$ is a countable support forcing iteration based on
$\langle Q_\eta\colon\allowbreak\eta<\kappa\rangle$.  Suppose for every $\eta<\kappa$ 
we have $P_\eta$ is weakly ${}^\omega\omega$-bounding and\/ {\rm
${\bf 1}\forces_{P_\eta}``Q_\eta$ is proper.''}
Then\/ $P_\kappa$ is weakly ${}^\omega\omega$-bounding.

Proof: Use the proof of Theorem 2.5 with $([\omega]^{\aleph_0},\allowbreak
{\cal P},\supseteq^*)$ replaced with
$({}^\omega\omega,{}^\omega\omega\cap V,\allowbreak\leq^*)$.

The Theorem is established.

The following definition is equivalent to [9, Definition VI.3.5(1)].

\proclaim Definition 3.5. Suppose $P$ is a forcing notion.  We say $P$ is almost 
${}^\omega\omega$-bounding iff\/ {\rm ${\bf 1}\forces``(\forall f\in{}^\omega\omega)\allowbreak(\exists g\in{}^\omega\omega\cap V)\allowbreak(\forall A\in[\omega]^{\aleph_0}\cap V)\allowbreak(\exists^\infty n\in A)\allowbreak(f(n)< g(n))$.''}

\proclaim Lemma 3.6.  Suppose $P$ is almost ${}^\omega\omega$-bounding.  Then $P$ is weakly ${}^\omega\omega$-bounding.

Proof: Take $A=\omega$ in Definition 3.5.

\proclaim Lemma 3.7.  Suppose $P$ is weakly ${}^\omega\omega$-bounding and\/ {\rm ${\bf 1}\forces_P``Q$ is almost ${}^\omega\omega$-bounding.''}  Then $P*Q$ is weakly ${}^\omega\omega$-bounding.

Proof: Like Lemma 2.7.

The Lemma is established.

\proclaim Theorem 3.9.  Suppose $\langle P_\eta\,\colon\allowbreak\eta\leq\kappa\rangle$ is
a countable support forcing iteration based on $\langle Q_\eta\,\colon\allowbreak\eta
<\kappa\rangle$.  Suppose for every $\eta<\kappa$ we have\/ {\rm 
${\bf 1}\forces_{P_\eta}``Q_\eta$
is proper and almost ${}^\omega\omega$-bounding.''}
Then $P_\kappa$ is weakly ${}^\omega\omega$-bounding.

Proof: By Theorem 3.4 and Lemma 3.7.

\section{Preservation of no new random reals}

We now turn our attention to the preservation of the property ``no new random reals.''

\proclaim Definition 4.1. For $\tau\in{}^{<\omega}2$, we let\/
 {\rm $U_\tau=\{\eta\in{}^\omega 2\,\colon\allowbreak\eta$ extends $\tau\}$.}

Recall that for $A\subseteq {}^\omega 2$, the outer measure of $A$ is\/ {\rm $\mu^*(A)={\rm inf}\{\Sigma\{2^{-{\rm lh}(\tau)}\,\colon\allowbreak \tau\in C\}\,\colon\allowbreak C\subseteq {}^{<\omega}2$ and $A\subseteq\bigcup\{U_\tau\,\colon\allowbreak\tau\in C\}\}$.}
$A$ is Lebesgue measurable iff $(\forall\tau\in{}^{<\omega}2)\allowbreak(\mu^*(A\cap U_\tau)+\mu^*(U_\tau-A)=\mu^*(U_\tau))$, in which case we write $\mu(A)=\mu^*(A)$.

\proclaim Definition 4.2.  Suppose $A\subseteq{}^\omega 2$.  We say that $A$ is closed under rational translation iff $(\forall b\in A)\allowbreak(\forall b^*=_{{\rm a.e.}}b)(b^*\in A)$.

The following Lemma is known as ``Kolmogorov's zero-one Law.''

\proclaim Lemma 4.3.  Suppose $A\subseteq{}^\omega 2$ is closed under rational translations and suppose that $A$ is Lebesgue measurable.  Then $\mu(A)=0$ or $\mu(A)=1$.

Proof:  Let $\gamma=\mu(A)$ and suppose, towards a contradiction, that $0<\gamma<1$.

Claim 1. Whenever $\tau\in{}^{<\omega}2$ and $\tau_0$ and $\tau_1$ are the immediate successors of $\tau$, then $\mu(A\cap U_{\tau_0})=\mu(A\cap U_{\tau_1})$.

Proof:  We have $\mu(A\cap U_{\tau_0})=2^{-{\rm lh}(\tau_0)}\mu(\{b\in {}^\omega 2\,\colon\allowbreak \tau_0\hat{\ }b\in A\})=2^{-{\rm lh}(\tau_1)}\mu(\{b\in {}^\omega 2\,\colon\allowbreak \tau_1\hat{\ }b\in A\})=\mu(A\cap U_{\tau_1})$.

Claim 2: For all $\tau\in{}^{<\omega}2$ we have $\mu(A\cap U_\tau)= 2^{-{\rm lh}(\tau)}\gamma$.

Proof:  By induction on $\tau$, using Claim 1.

Choose $\delta>\gamma$ such that $\delta^2<\gamma$.  Choose $C\subseteq{}^{<\omega}2$ such that $A\subseteq\bigcup\{U_\tau\,\colon\allowbreak\tau\in C\}$ and $\Sigma\{\mu(U_\tau)\,\colon
\allowbreak\tau\in C\}<\delta$.

For each $\tau\in C$, we may, using Claim 2,  choose $C_\tau\subseteq {}^{<\omega}2$ such that $A\cap U_\tau
\subseteq\bigcup\{U_\eta\,\colon\allowbreak\eta\in C_\tau\}$ and
$\Sigma\{\mu(U_\eta)\,\colon\allowbreak\eta\in C_\tau\}< 2^{-{\rm lh}(\tau)}\delta$.

Let $C^*=\bigcup\{C_\tau\,\colon\allowbreak\tau\in C\}$.

We have that $A\subseteq\bigcup\{U_\eta\,\colon\allowbreak\eta\in C^*\}$ and
$\Sigma\{\mu(U_\eta)\,\colon\allowbreak\eta\in C^*\}<\delta^2<\gamma$.

This contradiction establishes the Lemma.

\proclaim Definition 4.4.  Suppose $Y\subseteq{}^\omega 2$.  We define ${\rm RT}(Y)$, the ``rational translates'' of $Y$, to equal $\{b\in{}^\omega 2\,\colon\allowbreak
(\exists b'\in Y)\allowbreak(b'=_{\rm a.e.} b)\}$.

\proclaim Definition 4.5.  Suppose $y$ and $y'$ are perfect subsets of ${}^{\omega}2$
of positive Lebesgue measure.  We define $y\preceq y'$ to mean $y\subseteq
{\rm RT}(y')$.

\proclaim Lemma 4.6.  Suppose $\langle y_n\,\colon\allowbreak n\in\omega\rangle$ is a sequence of perfect subsets of ${}^\omega
2$ of positive Lebesgue measure.  Then there is a perfect set $y\subseteq{}^\omega 2$ of positive Lebesgue measure such
that $(\forall n\in\omega)\allowbreak(y\preceq y_n)$.

Proof: By Lemma 4.3 we have that $\mu({\rm RT}(y_n))=1$ for every $n\in\omega$.  For each $n\in\omega$ let $D_n
\subseteq{}^\omega 2$ be an open set such that $\mu(D_n)<2^{-n-1}$ and $D_n\cup {\rm RT}(y_n) = {}^\omega 2$.
Let $C={}^\omega 2 -\bigcup\{D_n\,\colon\allowbreak n\in\omega\}$. We have that $C$ is a closed set of positive measure.
Let $y$ be the perfect kernel
 of $C$ (see [5, page 66]).  We have that $y$ is a perfect set of positive measure, and for every
$n\in\omega$ we have $y\subseteq C\subseteq {}^\omega 2 - D_n\subseteq {\rm RT}(y_n)$.

The Lemma is established.

\proclaim Lemma 4.7.  Suppose $x$ and $y$ are subsets of ${}^\omega 2$. Then $x\cap{\rm RT}(y)=\emptyset$ iff
${\rm RT}(x)\cap y = \emptyset$.

Proof: Clear.

\proclaim Lemma 4.8. Let $P$ be any forcing. Then\/ {\rm $V[G_P]\models``(\forall
x\in{}^\omega 2)(x$ is random over $V$ iff
$(\forall y\in V)\allowbreak(y$ is a perfect set 
of positive Lebesgue measure implies $x\in{\rm RT}(y)))$.''}

Proof: Work in $V[G_P]$.  Suppose $x\in{}^\omega 2$ is not random over $V$.  Let
$B\in V$ be a Borel set such that $x\in B$ and $\mu(B)=0$. Let $D\in V$ be an open set such that $\mu(D)<1$ and
${\rm RT}(B)\subseteq D$.  Let $y$ be the perfect kernel of ${}^\omega 2-D$. Then $y\in V$ is
a perfect set of positive measure, and because $y\cap {\rm RT}(B) = \emptyset$, we have ${\rm RT}(y)\cap
B = \emptyset$, and therefore $x\notin {\rm RT}(y)$.

In the other direction, suppose $x\in{}^\omega 2$ and $y\in V$ is a perfect set of positive measure such that
$x\notin{\rm RT}(y)$.  We show that $x$ is not random over $V$.  Choose $\langle
D_n\,\colon\allowbreak n\in\omega\rangle\in V$ a sequence of open sets such that for every $n\in\omega$ we have
$\mu(D_n)<1/n$ and ${}^\omega 2-{\rm RT}(y)\subseteq D_n$.  Let $B=\bigcap\{D_n\,\colon\allowbreak
n\in\omega\}$.  We have that $B\in V$ is a Borel set of Lebesgue 
measure zero and $x\in B$.  Therefore $x$ is not random over $V$.

The Lemma is established.

The following is [9, Lemma VI.3.18].  Notice how the argument parallels the proof of Theorem 2.5.

\proclaim Theorem 4.9.  Suppose $\kappa $ is a limit ordinal and $\langle P_\eta\,\colon\allowbreak
\eta\leq\kappa\rangle $ is a countable support forcing iteration based on $\langle Q_\eta\,\colon\allowbreak
\eta<\kappa\rangle$. Suppose for every $\eta<\kappa$ we have\/ {\rm ${\bf 1}\forces_{P_\eta}``Q_\eta$
is proper and there are no reals that are random over $V$.''} Then\/ {\rm ${\bf 1}\forces_{P_\kappa}``$there
are no reals that are random over $V$.''}

Proof: For ${\rm cf}(\kappa) > \omega$ this is clear, so assume instead that ${\rm cf}(\kappa) =\omega$.

Suppose $p\in P_\kappa$ and $g$ is a $P_\kappa$-name and $p\forces``g\in{}^\omega 2$.'' Let $\lambda$ be
a sufficiently large regular cardinal and let $N$ be a countable elementary substructure of $H_\lambda$
containing $\{P_\kappa, p, g\}$.  Let $\langle\alpha_n\,\colon\allowbreak n\in\omega\rangle\in N$ be 
an increasing sequence cofinal in $\kappa$ such that $\alpha_0=0$.

Using Lemma 4.6, fix $y\subseteq{}^\omega 2$ a perfect set 
of positive Lebesgue measure such that for every perfect $y'\in N$
with $\mu(y')>0$ we have $y\preceq y'$.

Build $\langle q_k,p_k,m_k\,\colon\allowbreak k\in\omega\rangle$ such that $q_0=p$
 and for each
$k\in\omega$ we have the following:

(1) $p_k\in P_{\alpha_k}$ is $N$-generic, and

(2) $p_{k+1}\restr\alpha_k=p_k$, and

(3) $p_k\forces``q_{k+1}\in P_{\alpha_k,\kappa}\cap N[G_{P_{\alpha_k}}]$ and $q_{k+1}\leq q_k\restr[\alpha_k,\kappa)$,''
and

(4) $p_k\forces``p_{k+1}\restr[\alpha_k,\alpha_{k+1})\leq  q_{k+1}\restr\alpha_{k+1}$,'' and

(5) $p_k\forces``m_{k+1}>m_k$ and $q_{k+1}\forces``(\forall\rho\in{}^{m_k}2)\allowbreak
(U_{\rho\hat{\  }g\restr[m_k,m_{k+1})}\cap y=
\emptyset)$.'\thinspace''

The construction proceeds as follows.  Suppose we are given $p_k$ and $q_k$ and $m_k$.  Work in
$V[G_{P_{\alpha_k}}]$ with $p_k\in G_{P_{\alpha_k}}$. Build
$\langle q^i_k\,\colon\allowbreak
i\in\omega\rangle\in N[G_{P_{\alpha_k}}]$ a decreasing sequence of conditions in $P_{\alpha_k,\kappa}$ and
 $f_k\in{}^\omega 2\cap N[G_{P_{\alpha_k}}]$ such that $q^0_k\leq q_k$ and
for every $i\in\omega$ we have $q^i_k\forces``f_k(i)=g(i)$.'' Using the hypothesis on
$P_{\alpha_k}$ and Lemma 4.8, we may choose a perfect set $y_k\in V$ of positive measure
such that $f_k\notin{\rm RT}(y_k)$.  By elementarity we may assume $y_k\in N[G_{P_{\alpha_k}}]\cap
V = N$.

Because $y\preceq y_k$ we have
${\rm RT}(y)\subseteq{\rm RT}(y_k)$, and hence $f_k\notin {\rm RT}(y)$.
Hence by Lemma 4.7 we have ${\rm RT}(\{f_k\})\cap y = \emptyset$.
Hence for each $\rho\in {}^{m_k}2$, we may let $m_k^\rho$ be an integer greater
than $m_k$ such that $U_{\rho\hat{\  }f_k\restr[{\rm lh}
(\rho),m_k^\rho)}\cap y=\emptyset$, using the fact that $y$ is closed.

Let $m_{k+1} = {\rm max}\{m_k^\rho\,\colon\allowbreak\rho\in{}^{m_k} 2\}$. Let $q_{k+1}=
q_k^{m_{k+1}+1}$. We have that $q_{k+1}$ satisfies
(3) and (5).  Using the Proper Iteration Lemma, we may choose $p_{k+1}$ satisfying (1), (2), and (4).

This completes the recursive construction.

Let $r\in P_\kappa$ be chosen such that 
$(\forall k\in\omega)\allowbreak(r\restr\alpha_k=p_k)$.

We have $r\forces``{\rm RT}(\{g\})\cap y = \emptyset$.'' Hence by Lemmas 4.7 and 4.8 we have
$r\forces``g$ is not random over $V$.''

The Theorem is established.

\section{Preservation of ``every new dense open set contains an old dense open set''}

In this section we prove preservation of the property ``every new dense open set contains an old dense open set.'' Shelah includes two very different proofs of this fact in his book; we follow the proof given in [9, Section XVIII.3].

Throughout this section we fix an enumeration $\langle\eta_n^*\,\colon\allowbreak n\in\omega\rangle$ of ${}^{<\omega}\omega$ such that whenever $\eta^*_i$ is an initial segment of
$\eta^*_j$ then $i\leq j$.  Also, throughout this section we let ${\cal B}$ equal the set of functions from
${}^{<\omega}\omega$ into ${}^{<\omega}\omega$.  

\proclaim Definition 5.1.  Suppose $f$ and $g$ are in ${\cal B}$. 
We say $f\leq_{\cal B} g$ iff for every  $\eta\in{}^{<\omega}\omega$ there is
$\nu\in{}^{<\omega}\omega$ such that $\nu\hat{\  }f(\nu)$ is an initial segment of
$\eta\hat{\  }g(\eta)$.

We remark that Definition 5.1 differs from [9, Context and Definition XVIII.3.7A] because we have
incorporated
[9, Remark XVIII.3.7F(1)].  

\proclaim Lemma 5.2.  The relation $\leq_{\cal B}$ is a partial ordering of ${\cal B}$.

Proof: Immediate.

\proclaim Lemma 5.3.  Suppose $\langle f_i\,\colon\allowbreak i\in\omega\rangle$ is a sequence of elements of ${\cal B}$.  Then there is $g\in{\cal B}$ such that
for every $i\in\omega$ we have $f_{i}\leq_{\cal B}g$.

Proof: For every $\eta\in{}^{<\omega}\omega$ and $k\in\omega$ define
$g_0(\eta)=\eta$ and $g_{k+1}(\eta)=g_k(\eta)\hat{\  }f_k(\eta\hat{\  }g_k(\eta))$.
Define $g(\eta)$ to equal $g_n(\eta)$ where $\eta=\eta^*_n$.

To see that $f_k\leq_{\cal B}g$ it suffices to note that whenever $n>k$ then
\[(\eta^*_n\widehat{\  }g_k(\eta^*_n))\widehat{\  }f_k(\eta^*_n\widehat g_k(\eta^*_n))=
\eta^*_n\widehat{\  }g_{k+1}(\eta^*_n)\subseteq\eta^*_n\widehat{\  }g(\eta^*_n).\]
The Lemma is established.

\proclaim Lemma 5.4. Suppose $P$ is a forcing notion.  Then every dense open subset of\/ ${}^\omega\omega$ 
in $V[G_P]$ contains a dense open subset of\/ ${}^\omega\omega$ in $V$ iff\/ {\rm $V[G_P]\models``(\forall
f\in{\cal B})\allowbreak(\exists g\in{\cal B}\cap V)\allowbreak(
f\leq_{\cal B}g)$.''}

Proof: We first establish the ``if'' direction.  Work in $V[G_P]$. Suppose $D$ is a
dense open subset of ${}^\omega\omega$. Pick $f\in{\cal B}$ such that for every
$\eta\in{}^{<\omega}\omega$ we have $U_{\eta\hat{\  }f(\eta)}\subseteq D$.
Fix $g\in{\cal B}\cap V$ such that $f\leq_{\cal B}g$. 
Let $D'=
\bigcup\{U_{\eta\hat{\  }g(\eta)}\,\colon\allowbreak \eta\in{}^{<\omega}\omega\}$.
We have  $D'$ is a dense open subset of $D$ and $D'\in V$.

For the ``only if'' direction, suppose $f\in{\cal B}$.  
Build $\langle D_n, \eta_n, x_n\,\colon\allowbreak n\in\omega\rangle$ recursively such that for every
$n\in\omega$ we have that either $U_{\eta^*_n}\subseteq\bigcup\{D_i\,\colon i<n\}$ and
$D_n=D_{n-1}$ and $x_n=x_{n-1}$ and $\eta_n=\eta_{n-1}$, or all of the following::

(1) $\eta_n$ extends $\eta^*_n$, and

(2) $D_n = U_{\eta_n\hat{\  }f(\eta_n)\hat{\  }\langle 0\rangle}$, and

(3) $D_n$ is disjoint from $\bigcup\{D_i\,\colon i<n\}\cup\{x_i\,\colon\allowbreak i<n\}$, and

(4) 
$x_n\in {}^\omega\omega$ extends $\eta_n\hat{\  }f(\eta_n)\hat{\  }\langle 1\rangle$.

We may take $D'\in V$ open dense such that $D'\subseteq\bigcup\{D_n\,\colon\allowbreak n\in\omega\}$.

 Choose $g\in{\cal B}\cap V$ such that
$(\forall\eta\in{}^{<\omega}\omega)\allowbreak(U_{\eta\hat{\  }g(\eta)}\subseteq D')$.

Given $\eta\in{}^{<\omega}\omega$, pick $n\in\omega$  such that
$U_{\eta\hat{\  }g(\eta)}\cap U_{\eta_n\hat{\  }f(\eta_n)\hat{\  }\langle 0\rangle}\ne\emptyset$.
We have $x_n\notin D'$, and so $U_{\eta\hat{\  }g(\eta)}\subseteq U_{\eta_n\hat{\  }f(\eta_n)}$.
It follows that $f\leq_{\cal B}g$.

The Lemma is established.

The following is [9, Conclusion VI.2.15D] and [9, Claim XVIII.3.7D]; we follow the proof given in
[9, Chapter XVIII].

\proclaim Theorem 5.5.  Suppose  $\langle P_\eta\,\colon\eta\leq\kappa\rangle$ is a countable support forcing iteration based on
$\langle Q_\eta\colon\allowbreak\eta<\kappa\rangle$.  Suppose for every $\eta<\kappa$ 
we have\/ {\rm
${\bf 1}\forces_{P_\eta}``Q_\eta$ is proper and ${\bf 1}\forces_{Q_\eta}`$for every dense open $D\subseteq
{}^\omega\omega$ there is a dense open $D'\subseteq D$ such that $D'\in V[G_{P_\eta}]$.'\thinspace''}
Then\/ {\rm ${\bf 1}\forces_{P_\kappa}``$for every dense open $D\subseteq
{}^\omega\omega$ there is a dense open $D'\subseteq D$ such that $D'\in V$.''}

Proof: By induction on $\kappa$.
The induction step is clear for $\kappa$ a successor ordinal and, in light of Lemma
5.4, it is likewise clear for $\kappa$
 of uncountable cofinality. So we assume
$\cf(\kappa)=\omega$.

Suppose $p\in P_\kappa$ and $f$ is a $P_\kappa$-name and $p\forces``f\in{\cal B}$.'' Choose $\lambda$  a sufficiently large regular cardinal and
$N$  a countable elementary substructure of $H_\lambda$ such that $\{P_\kappa,f,p\}\in N$.

Let $\langle\alpha_k\,\colon\allowbreak k\in\omega\rangle\in N$ be an increasing sequence cofinal in $\kappa$
such that $\alpha_0=0$.
Using Lemma 5.3, 
fix $g\in{\cal B}$ such that $(\forall h\in{\cal B}\cap N)\allowbreak(h\leq_{\cal B} g)$.

Build $\langle q_k,p_k,m_k\,\colon k\in\omega\rangle$ such that $q_0=p$  and 
for every $k\in\omega$ we have that each of the following holds:

(1) $p_k\in P_{\alpha_k}$ is $N$-generic, and

(2) $p_k\forces``q_{k+1}\in P_{\alpha_k,\kappa}\cap N[G_{P_{\alpha_k}}]$ and $q_{k+1}\leq q_k\restr[\alpha_k,\kappa)$,'' and

(3) $p_k\forces``q_{k+1}\forces`m_k\in\omega$ and $\eta_k^*\hat{\  }g(\eta^*_k)$ extends $\eta^*_{m_k}\hat{\  }f(\eta^*_{m_k})$,'\thinspace'' and

(4) $p_{k+1}\restr\alpha_k =p_k$, and

(5) $p_k\forces``p_{k+1}\restr[\alpha_k,\alpha_{k+1})\leq q_{k+1}\restr\alpha_{k+1}$.''

The construction proceeds as follows. Given $p_k$ and $q_k$, work in 
$V[G_{P_{\alpha_k}}]$ with $p_k\in G_{P_{\alpha_k}}$.

Build $\langle q^i_k\,\colon\allowbreak i\in\omega\rangle\in N[G_{P_{\alpha_k}}]$ and
$f_k\in{\cal B}\cap N[G_{P_{\alpha_k}}]$ such that 
$q^0_k=q_{k}\restr[\alpha_k,\kappa)$, and for every $i\in\omega$ we have the following:

(1) $q^{i+1}_k\leq q^i_k$, and

(2) $q_k^{i+1}\forces``f_k(\eta^*_i)=f(\eta^*_i)$.''

Using Lemma 5.4, choose $g_k\in{\cal B}\cap V$ such that $f_k\leq_{\cal B}g_k$.
We may assume $g_k\in N[G_{P_{\alpha_k}}]$.  Hence $g_k\in N$.  Hence
$g_k\leq_{\cal B}g$.

By Lemma 5.2 we have $f_k\leq_{\cal B}g$, so we may 
 choose $m_k$  such that 
$\eta^*_k\hat{\  }g(\eta^*_k)$ extends $\eta^*_{m_k}\hat{\  }f_k(\eta^*_{m_k})$.

Let $q_{k+1}=q_k^{m_k+1}$.  We have that $q_{k+1}$ and $m_k$ satisfy (2) and (3).

Using the Proper Iteration Lemma we may choose $p_{k+1}$ satisfying (1), (4), and (5).

This completes the recursive construction.

Let $r\in P_\kappa$ be such that for every $k\in\omega$ we have $r\restr \alpha_k
=p_k$.

Suppose, towards a contradiction, that $r'\leq r$ and $r'\forces``f\not\leq_{\cal B}g$.''  
Fix a $P_\kappa$-name $k$ such that
$r'\forces``\eta^*_k\hat{\  }g(\eta^*_k)$ does not extend $\eta^*_{m_k}\hat{\  }f(\eta^*_{m_k})$.''

By strengthening $r'$ we may assume that $k$ and $m_k$ are integers rather than  names.

Because $r'\leq(p_{k+1},q_{k+1})$ we have
$r'\forces``\eta^*_k\hat{\  }g(\eta^*_k)$  extends $\eta^*_{m_k}\hat{\  }f(\eta^*_{m_k})$.''
This is a contradiction.

The Theorem is established.

\section{On ``the set of reals that are 
in the ground model has positive outer measure in the forcing extension''}

In this section we present a theorem of Shelah ([9, Claim XVIII.3.8B(3)]) 
that gives a sufficient
condition for a forcing iteration to satisfy $\mu^*({}^\omega 2\cap V)>0$.
This notion has been investigated also by [4].

\proclaim Definition 6.1.  We let ${\cal B}'$ be the set of functions $f$ from $\omega$ into
${}^{<\omega}2$ such that $\Sigma\{\mu(U_{f(m)})\,\colon\allowbreak m\in\omega\}\leq 1$.

\proclaim Lemma 6.2.  Suppose $g\in{}^\omega 2$ and\/ $\lambda $ is a sufficiently large
regular cardinal and $N$ is a countable elementary substructure of\/ $H_\lambda$.
Then $g $ is random over $N$ iff\/ {\rm $(\forall f\in{\cal B}'\cap N)\allowbreak
(\exists m\in\omega)\allowbreak(\forall i\geq m)\allowbreak(g$ does not extend $f(i))$.}

Proof: We first establish the ``only if'' direction.  Suppose $g\in
{}^\omega 2$ and $f\in{\cal B}'\cap N$ and $(\exists^\infty m\in\omega)\allowbreak
(g$ extends $f(m))$.
 Let $B=\{h\in{}^\omega 2\,\colon\allowbreak(\exists^\infty m\in\omega)\allowbreak
(f(m)$ is an initial segment of $h)\}$.
 Then $B\subseteq {}^\omega 2$ is a Borel set and $g\in B\in N$, and $\mu(B)=0$ because
for every $n\in\omega$ we have that $B$ is covered by $\bigcup\{U_{f(i)}\,\colon\allowbreak
i\geq n\}$, and ${\rm lim}_{n\rightarrow\infty}(\mu(\bigcup\{U_{f(i)}\,\colon\allowbreak
i\geq n\})=0$. Therefore $g$ is not random over $N$.

To prove the ``if'' direction, suppose that $g$ is not random over $N$.  We may choose $B\in N$ a Borel set of measure zero such that $g\in B\in N$. Let $\langle D_n\,\colon n\in\omega\rangle\in N$
be a sequence of open subsets of $ {}^\omega 2$  such that for every $n\in\omega$ we have $B\subseteq D_n$
and $\mu(D_n)<2^{-n}$. For each $n\in\omega$ choose $k_n\leq\omega$ and $\langle\eta^n_i\,\colon i<k_n\rangle$
a sequence of pairwise incomparable elements of ${}^{<\omega}2$ such that $D_n=\bigcup
\{U_{\eta^n_i}\,\colon\allowbreak i<k_n\}$.  Furthermore we may assume that
$\langle\langle\eta^n_i\,\colon\allowbreak i<k_n\rangle\,\colon\allowbreak n\in\omega\rangle$ is an
element of $N$.  Let $f\in N$ be a one-to-one function mapping $\omega$ onto
$\{\eta^n_i\,\colon\allowbreak i<k_n$ and $n\in\omega\}$.  Then we have that $f\in{\cal B}'$
and $(\exists^\infty m\in\omega)\allowbreak(g\in U_{f(m)})$.  The Lemma is established.

\proclaim Lemma 6.3.  Suppose $g\in{}^\omega 2$ and\/ $\lambda $ is a sufficiently large
regular cardinal and $N$ is a countable elementary substructure of\/ $H_\lambda$.
Suppose $g$ is random over $N$.  Suppose $Y\in N$ is a subset of\/ ${}^{<\omega}2$ and
$\Sigma\{\mu(U_{\eta})\,\colon\eta\in Y\}$ is finite.  Then $\{\eta\in Y\,\colon\allowbreak
g\in U_{\eta}\}$ is finite.

Proof:  We may assume $Y$ is infinite. Choose a finite integer $m$ and infinite sets
(not necessarily disjoint) $D_i\subseteq Y$ for $i<m$ such that
each $D_i$ is in $N$ and $\bigcup\{D_i\,\colon\allowbreak i<m\}=Y$ and for each $i<m$ we have
$\Sigma\{\mu(U_{\eta})\,\colon\allowbreak\eta\in D_i\}\leq 1$.
For each $i<m$ choose $f_i\in N$ such  that $f_i$ maps $\omega$ onto $D_i$.
By Lemma 6.2, for every $i<m$ there is $\beta_i\in\omega$ such that 
$(\forall j\geq \beta_i)\allowbreak(g$ does not extend $f_i(j))$.
 Hence $\{\eta\in Y\,\colon
\allowbreak g\in U_{\eta}\}\subseteq\bigcup\{\{f_i(j)\,\colon\allowbreak j<\beta_i\}\,\colon\allowbreak
i<m\}$, which is finite.  The Lemma is established.

\proclaim Lemma 6.4. Suppose $P$ is a poset such that whenever $\lambda$ is a sufficiently large
regular cardinal and $N$ is a countable elementary substructure of $H_\lambda$ and $P\in N$ and
$g\in{}^\omega 2$ and $g$ is random over $N$, then\/ {\rm $V[G_P]\models``g $ is random
over $N[G_P]$.''} Then\/ {\rm $V[G_P]\models``{}^\omega 2\cap V$ has positive outer measure.''}

Proof: Suppose, towards a contradiction, that in $V[G_P]$
we have that $B $ is a Borel 
subset of ${}^\omega 2$ such that ${}^\omega 2\cap V \subseteq B$ and $\mu(B)=0$.

In $V$, choose $\lambda$ a sufficiently large
regular cardinal and $N$ a countable elementary substructure of $H_\lambda$ such that $p\in N$ and
a name for $B$ is in $N$.  Let $g\in{}^\omega 2$ be
random over $N$.  By hypothesis, $V[G_P]\models``g$ is random over $N[G_P]$.''  Therefore
$V[G_P]\models``g\notin B$.'' This contradiction establishes the Lemma.

\proclaim Lemma 6.5. Suppose $P$ is a poset.  Suppose $\chi$ is a sufficiently large regular cardinal and\/
$\lambda$ is a regular cardinal sufficiently larger than $\chi$.  Suppose $N$ is a countable elementary substructure of\/
$H_\lambda$ and $N_1$ and $N_2$ are countable elementary substructures of $H_\chi$ and
$\chi\in N$ and $P\in N_1\in N_2\in N$.  Suppose also

\noindent (1) {\sl $G_1\subseteq P\cap N_1$ is an $N_1$-generic subset of\/ 
$P$, and }

\medskip

\noindent (2) {\sl  $p\in G_1$ and $G_1\in N$, and}

\medskip

\noindent(3) {\sl $\langle f_l\,\colon\allowbreak l\leq k\rangle\in N$ is a finite sequence
of $P$-names such that\/ {\rm $p\forces``f_l\in {\cal B}'\cap
N_1[G_P]$''} for all\/ $l\leq k$, and}

\medskip

\noindent (4) {\sl $g\in{}^\omega 2$ is random over $N$, and}

\medskip

\noindent (5) {\sl
$\langle\beta_l\,\colon\allowbreak l\leq k\rangle$ is a sequence of integers and 
 for all $l\leq k$ we have
$(\forall j\geq\beta_l)\allowbreak(g $ does not extend $f_l[G_1](j))$.}
That is, for every $j\geq\beta_l$ there is $p'\in G_1$ and $\rho\in{}^{<\omega} 2$
such that\/ $g$ does not extend  $\rho$  and\/ {\rm
$p'\forces``\rho=f_l(j)$.''}

\medskip

\noindent {\bf Then} {\sl there is $G_2\subseteq P\cap N_2$  an $N_2$-generic subset of\/ 
$P$ such that $p\in G_2$ and $G_2\in N$ and  for all $l\leq k$ we have\/ {\rm
$(\forall j\geq\beta_l)\allowbreak(g$ does not extend $f_l[G_2](j))$.}}

\medskip

Proof: Build $\langle p_n\,\colon\allowbreak n\in\omega\rangle\in N$ and $\langle m_n\,\colon\allowbreak
n\in\omega\rangle\in N$ and $\langle f^*_l\,\colon\allowbreak l\leq k\rangle\in N$ such that $p_0=p$ and
for each $n\in\omega$ we have each of the following:

(1) $p_n\in G_1$ and $p_{n+1}\leq p_n$, and

(2) $m_n$ is an integer such that $m_n\geq n$ and
$p_n\forces``\Sigma\{\mu(U_{f_l(i)})\,\colon\allowbreak i\geq m_n\}<2^{-n}$ for each $l\leq k$,'' and

(3) for every $l\leq k$ we have $f^*_l\in N$ maps $\omega$ into ${}^{<\omega}2$, and

(4) $p_n\forces``f_l\restr m_n = f^*_l\restr m_n$ for each $l\leq k$.''

Claim 1. For $l\leq k$ we have $f^*_l\in{\cal B}'$.

Proof.  Suppose, towards a contradiction, that $l\leq k$ and $m\in\omega$ and
$\Sigma\{\mu(U_{f^*_l(i)})\,\colon\allowbreak i<m\}>1$. Because
$p_m\forces``f_l\restr m=f^*_l\restr m$,'' it follows that
$p_m\forces``f_l\notin{\cal B}'$.'' This contradiction establishes the Claim.

Build $\langle p_{n,m}\,\colon\allowbreak m\in\omega, n\in\omega\rangle\in N$ and $\langle f^*_{l,n}\,\colon
\allowbreak l\leq k, n\in\omega\rangle\in N$ such that each of the following holds:

(1) for every $n\in\omega$ we have that $\langle p_{n,m}\,\colon\allowbreak m\in\omega\rangle$ is an $N_2$-generic
sequence for $P$ and $p_{n,0}=p_n$, and

(2) for every $l\leq k$ and $n\in\omega$ and $m\in\omega$ we have 
$p_{n,m}\forces``f^*_{l,n}\restr m=f_l\restr m$.''

Claim 2. For $l\leq k$ and $n\in\omega$ we have $f^*_{l,n}\in{\cal B}'$.

Proof: Similar to Claim 1.

Claim 3. For every $l\leq k$ and $n\in\omega$ we have
$\Sigma\{\mu(U_{f^*_{l,n}(i)})\,\colon\allowbreak i\geq m_n\}\leq 2^{-n}$.

Proof: Suppose $l$ and $n$ constitute a counterexample. Then we can choose an integer $t$ so large that
$\Sigma\{\mu(U_{f^*_{l,n}(i)})\,\colon\allowbreak m_n\leq i<t\}>2^{-n}$.
We have $p_{n,t}\forces``\Sigma\{\mu(U_{f^*_{l,n}(i)})\,\colon\allowbreak m_n\leq i<t\}=
\Sigma\{\mu(U_{f^*_{l}(i)})\,\colon\allowbreak m_n\leq i<t\}<
\Sigma\{\mu(U_{f^*_{l}(i)})\,\colon\allowbreak m_n\leq i<\omega\}<2^{-n}$.''
This contradiction establishes the Claim.

For each $l\leq k$ and $n\in\omega$ let $U^*_{l,n}=\bigcup\{U_{f^*_{l,n}(i)}\,\colon\allowbreak
i\in\omega\}$.

Claim 4. For every $l\leq k$ and $n\in\omega$ we have $U^*_{l,n}\subseteq
\bigcup\{U_{f^*_l(i)}\,\colon\allowbreak i\in\omega\}\cup\bigcup\{U_{f^*_{l,n}(i)}\,\colon\allowbreak
i\geq m_n\}$.

Proof: The Claim is forced by the condition $p_{n,n}$, hence it is true outright.

For each $l\leq k$ let $U^*_l=\bigcup\{U^*_{l,n}\,\colon\allowbreak n\in\omega\}$. By Claims
3 and 4 we have that $\mu(U^*_l)$ is finite for every $l\leq k$.
By Lemma 6.4 we have that
 $\{\rho\in{}^{<\omega}2\,\colon\allowbreak
(\exists l\leq k)\allowbreak(\exists n\in\omega)\allowbreak(\exists i\in\omega)\allowbreak
(\rho=f^*_{l,n}(i)$ and $g$ extends $\rho)\}$ is finite.
Therefore, we may fix
$n^*$ so large that $(\forall l\leq k)\allowbreak(\forall n\in\omega)\allowbreak
(\forall i\in\omega)(g$ extends $f_{l,n}^*(i)$ only if
$\mu(U_{f^*_{l,n}(i)})\geq 2^{-n^*})$.

Claim 5. Suppose  $l\leq k$ and $i\in\omega$ and
$n\in\omega$ and $\mu(U_{f^*_{l,n}(i)})\geq 2^{-n^*}$. Then $i<m_{n^*}$.

Proof: Suppose $i\geq m_{n^*}$. Then $p_{n,i+1}\forces``\mu(U_{f^*_{l,n}(i)})=
\mu(U_{f_l(i)})<\Sigma\{\mu(U_{f_l(j)})\,\colon\allowbreak j\geq m_{n^*}\}< 2^{-n^*}$.
This contradiction establishes the Claim.

 Fix $t>m_{n^*}$ such that $t>\beta_l$ for every $l\leq k$.
For every $l\leq k$ we have $p_{n^*,t}\forces``f^*_{l,n^*}\restr t=f^*_l\restr t$.''
Thus, by Claim 5,
 we have that $p_{n^*,t}\forces``(\forall l\leq k)\allowbreak(\forall i\geq\beta_l)\allowbreak
(g$ does not extend $f^*_{l,n^*}(i))$.''

Let $G_2=\{p'\in P\cap N_2\,\colon\allowbreak (\exists m\in\omega)\allowbreak(p_{n^*,m}\leq p')\}$.
We have that $G_2$ is as required.

The Lemma is established.

\proclaim Definition 6.6. Suppose\/ $g\in{}^\omega 2$. We say that\/
 $P$ is\/ $g$-good iff\/ {\bf whenever}

\noindent (1) {\sl $\chi$ is a sufficiently large regular
cardinal and\/ $\lambda$ is a regular cardinal sufficiently larger than\/ $\chi$ and}

\medskip

\noindent (2) {\sl
$N$ is a countable elementary substructure of\/ $H_\lambda$ and\/ $\chi\in N$ 
and\/  $N_1$ is a countable elementary  substructure
of\/ $H_\chi$
and}

\medskip

\noindent(3) {\sl $P\in N_1\in N$ and}

\medskip

\noindent(4) {\sl
$g$ is random over\/ $N$  and}

\medskip

\noindent (5) {\sl $k\in\omega$ and\/
$\langle f_l\,\colon\allowbreak l<k\rangle\in N$ is a sequence of\/ $P$-names and}

\medskip

\noindent (6) {\sl
 $p\in P\cap N_1$ and}

\medskip

\noindent(7) {\sl {\rm
$p\forces``(\forall l<k)\allowbreak(f_l\in{\cal B}'\cap N_1[G_P])$,''} and}

\medskip

\noindent(8) {\it $\langle f^*_l\,\colon\allowbreak
l<k\rangle$ is a sequence of elements of\/ ${\cal B}'$ and
 $\langle\beta_l\,\colon\allowbreak
l<k\rangle $ is a sequence of integers and for every $l<k$ we have\/ {\rm
$(\forall m\geq\beta_l)\allowbreak(g$ does not extend $f^*_l(m))$}
 and}

\medskip

\noindent(9) {\sl $G_1\subseteq P\cap N_1$
and
$G_1\in N$ and $G_1$ is $N_1$-generic over $P$ and $p\in G_1$ and}

\medskip

\noindent (10) {\sl
$(\forall l<k)\allowbreak(f_l[G_1]=f^*_l)$,}

\medskip

\noindent {\bf then} {\sl  there  is $q\leq p$ such that\/ $q$
is $N$-generic and\/ {\rm
$q\forces``g $ is random over $N[G_P]$ and $(\forall l<k)\allowbreak(\forall m\geq\beta_l)\allowbreak(
g$ does not extend $f_l(m))$.''}}

\proclaim Lemma 6.7.  Suppose we have that

\noindent(1) {\sl $g\in{}^\omega 2$ and}

\medskip

\noindent(2) {\sl  {\rm
${\bf 1}\forces``Q$ is $g$-good,''} and}

\medskip

\noindent (3) {\sl $\chi$ is a sufficiently
large regular cardinal and\/ $\lambda$ is a regular cardinal sufficiently larger
than $\chi$, and}

\medskip

\noindent (4) {\sl $N$ is a countable elementary substructure of\/ $H_\lambda$
 and\/
$\{P*Q,\chi\}\in N$, and}

\medskip

\noindent (5) {\sl $p\in P$ is $N$-generic and $q$ is a $P$-name 
and}

\medskip

\noindent (6) {\sl {\rm $p\forces``N_1$ is a countable elementary substructure
of $H_\chi[G_P]$ and $N_1\in N[G_P]$ and $g $ is random over $N[G_P]$ and $q\in Q\cap N_1$,''}
and}

\medskip

\noindent (7) {\sl $k\in\omega$ and\/ {\rm $p\forces``\langle f_l\,\colon\allowbreak l<k\rangle\in N[G_P]$ is a 
sequence of $Q$-names and $q\forces_{Q}`(\forall l<k)\allowbreak
(f_l\in{\cal B}'\cap N_1[G_{Q}])$.,'\thinspace''} and}

\medskip

\noindent (8) {\sl  $\langle f^*_l\,\colon\allowbreak l<k\rangle$ and
$\langle\beta_l\,\colon\allowbreak l<k\rangle$ are sequences of $P$-names and\/
{\rm $p\forces``(\forall l<k)\allowbreak(f^*_l\in{\cal B}'\cap N[G_P]$ and
$\beta_l\in\omega$ and $(\forall i\geq\beta_l)\allowbreak(g$ does not
extend $f^*_l(i)))$,''} and}

\medskip

\noindent (9) {{\sl  $G$ is a $P$-name and\/ {\rm $p\forces``G\subseteq
Q\cap N_1$ is generic over $N_1$ and $q\in G\in N[G_P]$ and
$(\forall l<k)\allowbreak(f^*_l=f_l[G])$.''}}

\medskip

\noindent {\bf Then} {\sl  there is a $P$-name $r$ such that\/ {\rm $p\forces``r\leq q$''}
and $(p,r)$ is $N$-generic and\/ {\rm $(p,r)\forces``g$ is
random over $N[G_{P*Q}]$ and $(\forall l<k)\allowbreak(\forall i\geq\beta_l)\allowbreak
(g$ does not extend $f_l(i))$.''}}

\medskip

Proof: Immediate.

\proclaim Theorem 6.8.  Suppose $g\in{}^\omega 2$ and suppose
 $\langle P_\eta\,\colon\eta\leq\kappa\rangle$ is a countable support forcing iteration based on
$\langle Q_\eta\colon\allowbreak\eta<\kappa\rangle$.  Suppose for every $\eta<\kappa$ 
we have\/ {\rm
${\bf 1}\forces_{P_\eta}``Q_\eta$ is proper and $g$-good.'\thinspace''}
Suppose also

\noindent (1) {\sl $\chi$ is a sufficiently large regular cardinal and\/
 $\lambda$ is a regular cardinal sufficiently larger than $\chi$, and}

\medskip

\noindent(2) {\sl  $N$ is a countable elementary substructure of\/ $H_\lambda$ and $\{P_\kappa,\chi\}\in N$
and}

\medskip

\noindent (3) {\sl   $\alpha\in \kappa\cap N$ and $p\in P_\alpha\cap N$
and}

\medskip

\noindent (4) {\sl  {\rm $p\forces``N'$ is a countable
elementary substructure of $H_\chi[G_{P_\alpha}]$ and $P_{\alpha,\kappa}\in N'\in N[G_{P_\alpha}]$
(so necessarily $\alpha\in N'$),''} and}

\medskip

\noindent (5) {\sl 
 $p$ is $N$-generic and $g$ is a $P_\alpha$-name and\/ 
{\rm $p\forces``g$
is random over $N[G_{P_\alpha}]$ and $q\in P_{\alpha,\kappa}\cap N'$,''} and}

\medskip

\noindent (6) {\sl
 $k\in\omega$ and\/ {\rm $p\forces``\langle f_l\,\colon\allowbreak l<k\rangle\in N[G_{P_\alpha}]$ is a 
sequence of $P_{\alpha,\kappa}$-names and $q\forces_{P_{\alpha,\kappa}}`(\forall l<k)\allowbreak
(f_l\in{\cal B}'\cap N'[G_{P_{\alpha,\kappa}}])$,'\thinspace''} and}

\medskip

\noindent (7) {\sl
 $\langle f^*_l\,\colon\allowbreak l<k\rangle$ and
$\langle\beta_l\,\colon\allowbreak l<k\rangle$ are sequences of $P_\alpha$-names and\/
{\rm $p\forces``(\forall l<k)\allowbreak(f^*_l\in{\cal B}'\cap N[G_{P_\alpha}]$ and
$\beta_l\in\omega$ and $(\forall i\geq\beta_l)\allowbreak(g$ does not
extend $f^*_l(i)))$,''} and}

\medskip

\noindent (8) {\sl   $G$ is a $P_\alpha$-name and\/ {\rm $p\forces``G\subseteq
P_{\alpha,\kappa}\cap N'$ is generic over $N'$ and 
$q\in G\in N[G_{P_\alpha}]$ and
$(\forall l<k)\allowbreak(f^*_l=f_l[G])$.''}

\medskip

\noindent {\bf
Then} {\sl  there is $r\in P_\kappa$ such that\/ $r\restr\alpha=p$ and\/ 
{\rm $p\forces``r\restr[\alpha,\kappa)\leq q$''}
and $r$ is $N$-generic and\/ {\rm $r\forces``g$ is
random over $N[G_{P_\kappa}]$ and $(\forall l<k)\allowbreak(\forall i\geq\beta_l)\allowbreak
(g$ does not extend $f_l(i))$.''}}}

\medskip

Proof: By induction on $\kappa$.

Successor case: $\kappa = \gamma+1$.

In $V[G_{P_\alpha}]$ let $G_1=G\restr\gamma$ and $G_2=G/G_1$.  That is, $G_1=\{p'\restr\gamma\,\colon\allowbreak
p'\in G\}$ and $(\forall p'\in P_{\alpha,\gamma})\allowbreak(\forall r')\allowbreak(p'\forces``
r'\in G_2$'' iff $(\forall p^*\leq p')\allowbreak(\exists q'\in G)\allowbreak(\exists p^\#\leq p^*)
\allowbreak(p^\#\leq q'\restr\gamma$ and $p^\#\forces``r'=q'(\gamma)$'')).

Choose $N^*$ a countable elementary substructure of $H_\chi[G_{P_\alpha}]$ such that
$N'[G_{P_\alpha}]\in N^*\in N[G_{P_\alpha}]$ and  $G\in N^*$.
Choose $\langle f^{**}_l\,\colon\allowbreak l<k\rangle$ such that for
all $l<k$ we have $(p,q\restr\gamma)\forces_{P_\gamma}``f_l^{**}=f_l[G_2]$.''
Because $p\forces``f^{**}_l[G_1]=f^*_l$'' for all $l<k$, we have that
$(p,q\restr\gamma)\forces``(\forall l<k)\allowbreak(\forall j\geq
\beta_l)\allowbreak(g$ does not extend $f^{**}_l(j))$.''

Use Lemma 6.5 to choose $G_1'$ such that $p\forces``G_1'\subseteq 
P_{\alpha,\gamma}\cap N^*$ is generic over $N^*$  and
$q\restr\gamma\in G'_1$ and $G'_1\in N[G_{P_\alpha}]$ and $(\forall l<k)\allowbreak
(\forall j\geq\beta_l)\allowbreak(g$ does not extend $f_l^{**}[G'_1](j))$.''

By the induction hypothesis, with $G_1'$ playing the role of $G$ and $\langle
f^{**}_l\,\colon l<k\rangle$ playing the role of $\langle f_l\,\colon\allowbreak l<k\rangle$,
 we can choose $r'\in P_\gamma$
such that
$r'\restr\alpha=p$ and\/ 
{\rm $p\forces``r'\restr[\alpha,\gamma)\leq q\restr\gamma$''}
and $r'$ is $N$-generic and\/ {\rm $r'\forces``g$ is
random over $N[G_{P_\gamma}]$ and $(\forall l<k)\allowbreak(\forall i\geq\beta_l)\allowbreak
(g$ does not extend $f_l^{**}(i))$.''}

Using Lemma 6.7 with $G_2$ playing the role of $G$ and $N'[G_{P_\gamma}]$ playing
the role of $N_1$, we may choose $r^*$ such that 
$r'\forces``r^*\in  
Q_\gamma$ and $r^*\leq q(\gamma)$'' and 
 $(r',r^*)\forces``g$ is
random over $N[G_{P_\kappa}]$ and $(\forall l<k)\allowbreak(\forall i\geq\beta_l)\allowbreak
(g$ does not extend $f_l(i))$.''

Let $r=(r',r^*)$. This concludes the verification of the successor case.

Limit case: $\kappa$ is a limit ordinal.

Let $\langle \alpha_n\,\colon\allowbreak n\in\omega\rangle$ be an increasing sequence from $\kappa\cap N$
cofinal in ${\rm sup}(\kappa\cap N)$ such that $\alpha_0=\alpha$. Let $\langle \sigma_n\,\colon\allowbreak
n\in\omega\rangle$ list all $P_\kappa$-names $\sigma$ such that $\sigma\in N$ and
${\bf 1}\forces_{P_\kappa}``\sigma$ is an ordinal.'' Let $\langle f_l\,\colon\allowbreak
l\in\omega\rangle $ be a sequence that extends $\langle f_l\,\colon\allowbreak l< k\rangle$,
such that it lists the set of all $P_\kappa$-names $f$ in $N$ such that $(p,q)\forces``f\in {\cal B}'$.''

Build $\langle p_n,q_n,\beta_n,G_n,G^*_n,G'_n,N_n\rangle$ 
such that $p_0=p$ and $q_0=q$ and $G_0=G$ and $N_0=N'$ and
$\langle\beta_l\,\colon l\in\omega\rangle$ extends $\langle \beta_l\,\colon\allowbreak l<k\rangle$,
and for every $n\in\omega$ we have that each of the following holds:

(1) $p_n\forces``G'_{n}=G_n\restr\alpha_{n+1}$ and $G_{n+1}=G_n/G'_{n}$ 
(see the successor case, above),'' and

(2) $p_n\forces``N_{n+1}$ is a countable elementary substructure of
$H_\chi[G_{P_{\alpha_n}}]$ and $\{N_n[G_{P_{\alpha_{n-1},\alpha_{n}}}],\allowbreak
G_n, f_n, \alpha_{n+1},\sigma_n\}\in N_{n+1}\in N[G_{P_{\alpha_n}}]$'' (if $n=0$ then
replace $N_0[G_{P_{\alpha_{-1},\alpha_0}}]$ with $N_0$), and

(3) $\beta_n$ is a $P_{\alpha_n}$-name for an integer and
$p_n\forces``(\forall j\geq\beta_n)\allowbreak (g$ does not extend
$f_n[G_{n}](j))$,'' and

(4) $p_n\forces``G^{*}_n\subseteq P_{\alpha_n,\alpha_{n+1}}\cap N_{n+1}$
is $N_{n+1}$-generic and $q_n\restr\alpha_{n+1}\in G^*_{n+1}\in N[G_{P_{\alpha_n}}]$
and $(\forall l<{\rm max}(n+1,k))\allowbreak
(\forall j\geq\beta_l)\allowbreak(g$ does not extend $f_l[G_{n+1}][G^*_n](j))$,'' and

(5) $p_{n+1}\in P_{\alpha_{n+1}}$ is $N$-generic and
$p_{n+1}\forces``g$ is random over $N[G_{P_{\alpha_{n+1}}}]$ and
$(\forall i<{\rm max}(n+1,k))\allowbreak(\forall j\geq\beta_l)\allowbreak
(g$ does not extend $f_l[G_{n+1}](j))$,'' and

(6) $p_n\forces``p_{n+1}\restr[\alpha_n,\alpha_{n+1})\leq
q_n\restr\alpha_{n+1}$,'' and

(7) $p_{n+1}\forces``q_{n+1}\leq q_n\restr[\alpha_{n+1},\kappa)$ and
$q_{n+1}\in G_{n+1}$ and $q_{n+1}$ decides the value of $\sigma_n$ and
$q_{n+1}$ decides the value of $f_l\restr n$ for every $l\leq n$.''

The construction proceeds as follows. Given $p_n$ and $q_n$ and $G_n$, construct $G_n'$ and
$G_{n+1}$ as in (1) (see successor case, above). There is no problem in
choosing $N_{n+1}$ as in (2). We have that
$p_n\forces``f_n[G_{n}]\in{\cal B}'$'' by the reasoning of
Claim 1 in the proof of Lemma 6.5, hence
we may choose $\beta_n$ as in (3) because of Lemma 6.2.
We may choose $G^*_n$ as in (4) by Lemma 6.5.  We may choose $p_{n+1}$ satisfying
(5) and (6) by using the induction hypothesis.  There is no difficulty in choosing
$q_{n+1}$ satisfying (7).

Take $r\in P_\kappa$ such that for every $n\in\omega$ we have $r\restr\alpha_n=p_n$.

Claim. $r\forces``g$ is random over $N[G_{P_\kappa}]$.

Proof: Suppose not.  By Lemma 6.2 we may take
$r'\leq r$ and $l\in\omega$ such that $r'\forces``(\exists^\infty m\in\omega)\allowbreak(
g$ extends $f_l(m))$.''  By strengthening $r'$ further, we may assume there is an integer
$\beta^*$ such that $r'\forces``\beta_l=\beta^*$.'' 
By a further strengthening of $r'$ we may assume there is an integer $j\geq \beta^*$ such that
$r'\forces``g$ extends $f_l(j)$.'' Let $n={\rm max}(j+1,l+1)$.
By (7) we have that $p_{n+1}\forces``q_{n+1}\forces`f_l[G_{n+1}](j)=f_l(j)$.'\thinspace''
We have $p_{n+1}\forces``g$ does not extend $f_l[G_{n+1}](j)$.''
The Claim is established.

We have that $r$ is $N$-generic by the usual argument on ordinal names in $N$, and
it is clear that $r\forces``(\forall l<k)\allowbreak(\forall j\geq\beta_l)\allowbreak
(g$ does not extend $f_l(j))$.''

The Theorem is established.

The following Theorem is [9, Claim XVIII.3.8C(1)].

\proclaim Theorem 6.9.  Suppose $\langle P_\eta\,\colon\allowbreak\eta\leq\kappa\rangle$ is
a countable support iteration based on $\langle Q_\eta\,\colon\allowbreak
\eta<\kappa\rangle$ and for every $\eta<\kappa$ we have\/ {\rm ${\bf 1}\forces_{P_\eta}``Q_\eta$
is proper and for every $g\in{}^\omega 2$ we have that $Q_\eta$ is $g$-good.''}
Then\/ {\rm $V[G_{P_\kappa}]\models``{}^\omega 2\cap V$ does not have measure zero.''}

Proof: By Theorem 6.8 with $\alpha=k=0$ and Lemma 6.2.

\section{Preservation of ``the set of old reals is non-meager''}

Let ${\cal B}^*$ be the set of functions from ${}^{<\omega} 2$ into ${}^{<\omega} 2$.

\proclaim Definition 7.1.  Suppose $f\in{\cal B}^*$ and $g\in{}^\omega 2$. 
We say $f R^\dag g$ iff\/ {\rm $(\exists^\infty m\in\omega)\allowbreak
(g\restr m\hat{\  }f(g\restr m)$ is
an initial segment of $g)$.}

\proclaim Lemma 7.2. Suppose $X\subseteq{}^\omega 2$. Then $X$ is non-meager iff for every
$f\in{\cal B}^*$ there is $g\in X$ such that $f R^\dag g$.

Proof: Suppose $X$ is non-meager, and suppose $f\in{\cal B}'$.

For every $i\in\omega$ let $D_i=\bigcup\{U_{\tau\hat{\  }f(\tau)}\,\colon
\allowbreak(\exists n > i)\allowbreak(\tau\in{}^n 2)\}$.
 We have that each $D_i$ is an open dense set, so because $X$ is non-meager, we may
fix $g\in X\cap\bigcap\{D_i\,\colon\allowbreak i\in\omega\}$.
Clearly $f R^\dag g$.

For the converse, suppose $(\forall f\in{\cal B}^*)\allowbreak(\exists g\in X)\allowbreak
(f R^\dag g)$, and suppose $\langle D_i\,\colon\allowbreak i\in\omega\rangle$ is a
decreasing sequence of open dense subsets of ${}^\omega 2$.  We show $X\cap\bigcap\{D_i\,\colon\allowbreak
i\in\omega\}$ is non-empty. It suffices to find $g\in X$ such that
$(\exists^\infty j\in\omega)\allowbreak(g\in D_j)$.

Choose $f\in{\cal B}^*$ such that for every $\eta\in{}^{<\omega} 2$ we have
$U_{\eta\hat{\  }f(\eta)}\subseteq D_{{\rm lh}(\eta)}$.
Fix $g\in X$ such that $f R^\dag g$. Given $i\in\omega$ choose
$j>i$ such that
$g\restr j\hat{\  }f(g\restr j)$ is an initial segment of $g$.
Let $\eta=g\restr j$.
Then $g\in U_{\eta\hat{\  }f(\eta)}\subseteq D_j$.

The Lemma is established.

\proclaim Lemma 7.3. Suppose $\lambda$ is a sufficiently large regular cardinal and
$N$ is a countable elementary substructure of $H_\lambda$, and suppose $g\in{}^\omega 2$.
The following are equivalent:

(1) $(\forall f\in{\cal B}^*\cap N)(f R^\dag g)$.

(2) $(\forall f\in{\cal B}^*\cap N)(\exists m\in\omega)(g\restr m\hat{\  }f(g\restr m)$ is
an initial segment of $g)$.

(3) $g$ is Cohen over $N$.

Proof: It is obvious that (1) implies (2).

Suppose (2) holds and $D\in N$ is an open dense subset of ${}^\omega 2$. 
Choose $f\in {\cal B}^*\cap N$ such that $(\forall\nu\in{}^{<\omega}2)\allowbreak
(U_{\nu\hat{\  }f(\nu)}\subseteq D)$.  Using (2), choose $m\in\omega$ such that
$g\restr m\hat{\  }f(g\restr m)$ is an initial segment of $g$.  We have $g\in U_{g\restr m\hat{\  }
f(g\restr m)}\subseteq D$.  We conclude that $g\in\bigcap\{D\in N\,\colon\allowbreak D $ is an
open dense subset of ${}^\omega 2\}$, i.e., $g$ is Cohen over $N$.

Finally, suppose (3) holds and $f\in{\cal B}^*\cap N$.  Suppose $k\in\omega$.  Let $D_k=\{h\in{}^\omega 2
\,\colon\allowbreak(\exists m>k)\allowbreak(h\restr m\hat{\  }f(h\restr m)$ is an initial segment of $h)\}$.
It is easy to see that for every $k\in\omega$ we have $D_k$ is an open dense subset of ${}^\omega 2$.
Because $(\forall k\in\omega)\allowbreak(g\in D_k)$ we have that $f R^\dag g$.

The Lemma is established.

The following Lemma, due to Goldstern and Shelah, is [9, Lemma XVIII.3.11].

\proclaim Lemma 7.4. Suppose $P$ is a Suslin proper forcing (see {\rm [1, Section 7]}) and for 
every forcing $Q$ we have\/ {\rm
${\bf 1}\forces_Q``P$ is Suslin proper and ${\bf 1}\forces_P`{}^\omega 2\cap V[G_Q]$ is not meager.'\thinspace''}
Suppose $\lambda$ is a sufficiently large regular cardinal and $N$ is a countable elementary submodel
of $H_\lambda$ and $P\in N$ and $p\in P\cap N$ and $g\in{}^\omega 2$ is Cohen over $N$.  Then there is
$q\leq  p$ such that $q$ is $N$-generic and\/ {\rm $q\forces``g$ is Cohen over $N[G_P]$.''}

The proof presented in [9] is quite clear, so we do not repeat it here.

\proclaim Lemma 7.5.  Suppose
 $\langle P_\eta\,\colon\eta\leq\kappa\rangle$ is a countable
 support forcing iteration based on
$\langle Q_\eta\colon\allowbreak\eta<\kappa\rangle$.  
Suppose for every $\eta<\kappa$ 
we have\/ {\rm
${\bf 1}\forces_{P_\eta}``Q_\eta$ is a Suslin proper forcing and for every forcing $Q$ we have 
${\bf 1}\forces_Q`Q_\eta$ is Suslin proper and 
${\bf 1}\forces_{Q_\eta}``{}^\omega 2\cap V[G_{P_\eta}][G_Q]$ 
is not meager.''\thinspace'\thinspace''} Suppose $\lambda$ is a sufficiently large regular cardinal and
$N$ is a countable elementary substructure of $H_\lambda$ and $P_\kappa\in N$ and $\alpha\in
\kappa\cap N$ and $p\in P_\alpha$ is $N$-generic and\/ {\rm $p\forces``q\in P_{\alpha,\kappa}\cap
N[G_{P_\alpha}]$ and $g\in {}^\omega 2$ is Cohen over $N[G_{P_\alpha}]$.''}
Then there is $r\in P_\kappa$ such that $r$ is $N$-generic and $r\restr \alpha=p$ and $p\forces``r\restr[\alpha,\kappa)\leq
q$'' and\/ {\rm $r\forces``g$ is Cohen over $N[G_{P_\kappa}]$.''}

Proof: By induction on $\kappa$.

Case 1: $\kappa$ is a successor ordinal.

Let $\beta$ be the immediate predecessor of $\kappa$.  By the induction hypothesis we may take
$r'\in P_\beta$ such that $r'$ is $N$-generic and $r'\restr\alpha=p$
and $p\forces``r'\leq q\restr\beta$'' and $r'\forces``g$ is Cohen over $N[G_{P_\beta}]$.''
By Lemma 7.4 we may take $r^*\in Q_\beta$ such that $r'\forces``r^*\leq
q(\beta)$ and $r^*$ is $N[G_{P_\beta}]$-generic and $r^*\forces`g$ is Cohen over $N[G_{P_\beta}][Q_\beta]
$.'\thinspace''  Let $r\in P_\kappa$ be defined by $r\restr\beta=r'$ and
$r(\beta)=r^*$.  We have that $r$ satisfies the requirements of the Lemma.

Case 2: $\kappa$ is a limit ordinal.

Let $\langle\alpha_k\,\colon\allowbreak k\in\omega\rangle$ be an increasing sequence
from $\kappa\cap N$ cofinal in ${\rm sup}(\kappa\cap N)$
such that $\alpha_0=\alpha$.
Let $\langle \sigma_k\,\colon\allowbreak k\in\omega\rangle$ list all $P_\kappa$-names $\sigma$ in $N$ such
that ${\bf 1}\forces_{P_\kappa}``\sigma$ is an ordinal.''

Let $\langle f_i\,\colon\allowbreak i\in\omega\rangle $ list all $P_\kappa$-names $f$ in $N$
such that $V[G_{P_\kappa}]\models``f \in {\cal B}^*$,'' and let $\langle\eta'_m\,\colon\allowbreak
m\in\omega\rangle$ list ${}^{<\omega}2$.

Build $\langle q_k,p_k,n_k\,\colon k\in\omega\rangle$ such that $p_0=p$ and $q_0=q$ and 
for every $k\in\omega$ we have that each of the following holds:

(1) $p_k\in P_{\alpha_k}$ is $N$-generic, and

(2) $p_k\forces``q_{k+1}\in P_{\alpha_k,\kappa}\cap N[G_{P_{\alpha_k}}]$ and $q_{k+1}\leq q_k\restr[\alpha_k,\kappa)$,'' and

(3) $p_k\forces``q_{k+1}\forces`g$ is Cohen over $N[G_{P_{\alpha_k}}]$ and
$\sigma_k\in N$ and $n_k\in\omega$ and
$g\restr n_{k}\hat{\  }f_k(g\restr n_{k})$ is an initial segment
of $g$,'\thinspace'' and

(4) $p_{k+1}\restr\alpha_k =p_k$, and

(5) $p_k\forces``p_{k+1}\restr[\alpha_k,\alpha_{k+1})\leq q_{k+1}\restr\alpha_{k+1}$.''

The construction proceeds as follows. Given $p_k$ and $q_k$,  work in 
$V[G_{P_{\alpha_k}}]$ with $p_k\in G_{P_{\alpha_k}}$.

Build $\langle q_k^m\,\colon\allowbreak
m\in\omega\rangle\in N[G_{P_{\alpha_k}}]$ and $f'_k\in{\cal B}^*\cap N[G_{P_{\alpha_k}}]$
such that $\langle q_k^m\,\colon\allowbreak
m\in\omega\rangle$ is
a decreasing sequence of elements of $P_{\alpha_k,\kappa}$
and $q_k^0\leq q_k\restr[\alpha_k,\kappa)$ and there is an ordinal
$\tau$ such that $q^0_k\forces``\tau=\sigma_k$,'' and for every $m\in\omega$ we have that
$q_k^m\forces``f'_k(\eta'_m)=f_k(\eta'_m)$.'' Necessarily $\tau\in N[G_{P_{\alpha_k}}]$ and
therefore, because $p_k\in G_{P_{\alpha_k}}$ is $N$-generic, we have $\tau\in N$.
Because $g$ is Cohen over $N[G_{P_{\alpha_k}}]$ we may use Lemma 7.3 to take $n_k$
such that $g\restr n_k\hat{\  }f_k'(g\restr n_k)$ is an initial segment of $g$.

Let $q_{k+1}=q_k^{n_k+1}$.

Using the induction hypothesis, we may choose $p_{k+1}$ as required.

This completes the recursive construction.

Let $r\in P_\kappa$ be such that for every $k\in\omega$ we have $r\restr \alpha_k
=p_k$.

We have that $r$ is $N$-generic, because for each $k\in\omega$ we have $p_{k+1}\forces``q_{k+1}\forces`
\sigma_k\in N$.'\thinspace''

Suppose, towards a contradiction, that $r'\leq r$ and $r'\forces``g$ is
not Cohen over $N[G_{P_\kappa}]$.''  Choose $r^*\leq r'$ and $k\in\omega$ such that
$r\forces``(\forall m\in\omega)\allowbreak(g\restr m\hat{\  }f_k(g\restr m)$ is not an initial
segment of $g)$.''

Because $r^*\leq(p_{k+1},q_{k+1})$ we have
$r^*\forces``
g_{n_k}\hat{\  }f_k(g\restr n_k)$ is an initial segment of $g$.''
This is a contradiction.

The Lemma is established.

The following Theorem is [9, Claim XVIII.Claim 3.10C].

\proclaim Theorem 7.6.  Suppose
 $\langle P_\eta\,\colon\eta\leq\kappa\rangle$ is a countable
 support forcing iteration based on
$\langle Q_\eta\colon\allowbreak\eta<\kappa\rangle$.  Suppose for every $\eta<\kappa$ 
we have\/ {\rm
${\bf 1}\forces_{P_\eta}``Q_\eta$ is Suslin proper forcing and for every forcing $Q$ we have 
${\bf 1}\forces_Q`Q_\eta$ is Suslin proper and 
${\bf 1}\forces_{Q_\eta}``{}^\omega 2\cap V[G_{P_\eta}][G_Q]$ is not meager.''\thinspace'\thinspace''} Then\/ {\rm ${\bf 1}\forces_{P_\kappa}``{}^\omega 2\cap V$ is not
meager.''}

Proof:  Suppose, towards a contradiction, that $q\in P_\kappa$ and $q\forces``{}^\omega 2
\cap V$ is meager.'' By Lemma 7.3 we may take $f$ a $P_\kappa$-name for
an element of ${\cal B}^*$ such that
 $q\forces``(\forall g\in {}^\omega 2\cap V)\allowbreak(f R^\dag g$ fails).''
Let $\lambda$ be a sufficiently large regular cardinal and let $N$ be a countable elementary substructure
of $H_\lambda$ such that $\{P_\kappa,q,f\}\in N$.

Let $g\in{}^\omega 2$ be Cohen over $N$.

By Lemma 7.5 with $\alpha=0$ we may take $r\leq q$ such that $r\forces``g$ is Cohen over $N[G_{P_\kappa}]$.''
By Lemma 7.3 we have $r\forces``f R^\dag g$.''  This contradiction establishes the Theorem.

\section{References}

\medskip

[1]  Bartoszynski, T., and H. Judah, Set Theory: On the Structure of the Real Line, A K Peters, 1995.

\medskip

[2]  Goldstern, M., Tools for Your Forcing Construction, Set theory of the reals (Haim Judah, editor), Israel Mathematical Conference Proceedings, vol. 6, American Mathematical Society, pp. 305--360. (1993)

\medskip

[3]  Goldstern, M., and J. Kellner, New reals: Can live with them, can live without them, Math. Log. Quart. 52, No. 2, pp.~115--124, 2006.

 \medskip

[4]  Kellner, J., and S. Shelah, Preserving preservation, Journal of Symbolic Logic,
vol.~{\bf 70}, pp.~914--945, 2005.

\medskip

[5] Moschovakis, Y., Descriptive Set Theory, Studies in Logic and the Foundations of Mathematics, v.~100, North-Holland, 1980.

\medskip

[6] Schlindwein, C., Shelah's work on non-semi-proper iterations, I, Archive for Mathematical Logic, vol.~{\bf 47}, pp.~579--606, 2008.

\medskip

[7] Schlindwein, C., Shelah's work on non-semi-proper iterations, II, Journal of Symbolic Logic, vol.~{\bf 66}, pp.~1865--1883, 2001.

\medskip

[8] Schlindwein, C., Understanding preservation theorems: Chapter VI of {\em Proper and Improper Forcing},''
submitted.

\medskip

[9] Shelah, S., {\bf Proper and Improper Forcing}, Perspectives in Mathematical Logic, Springer, Berlin, 1998.

\end{document}